\documentclass[preprint,12pt]{elsarticle}
\usepackage{amssymb}
\usepackage{lineno}
\usepackage[tbtags]{amsmath}
\usepackage{graphicx,floatrow,threeparttable}

\newtheorem{lemma}{\indent Lemma}

\begin{document}
\begin{frontmatter}
\title{Equivalence of the Initialized Riemann-Liouville Derivatives and the Initialized Caputo Derivatives}
\author[rvt]{Jian Yuan \corref{cor}}
\address{College of Mathematic and Information Science, Shandong Technology and Business University, Yantai 264005, P.R.China}
\cortext[cor]{Corresponding author. Email: yuanjianscar@gmail.com}
\author{Guozhong Xiu, Bao Shi}
\address{Institute of System Science and Mathematics, Naval Aeronautical University, Yantai 264001, P.R.China}

\begin{abstract}
Initialization of fractional differential equations remains an ongoing problem. In recent years, the initialization function approach and the infinite state approach provide two effective ways to deal with this problem. The purpose of this paper is to prove the equivalence of the initialized Riemann-Liouville derivatives and the initialized Caputo derivatives with arbitrary orders. By synthesizing the above two initialization theories, the diffusive representations of the two initialized derivatives with arbitrary orders are derived. Laplace transforms of the two initialized derivatives are shown to be equal. As a result, the two most commonly used derivatives are proved to be equivalent when initial conditions are properly imposed.
\end{abstract}


\end{frontmatter}

\section{Introduction}
Fractional calculus provides a powerful tool to model real-world phenomena exhibiting memory and hereditary properties [1]. Initial conditions of fractional differential equations characterize the past histories of dynamical systems from the starting time to the initial time. As a result, initial conditions of fractional differential equations should be imposed in a different way from that of integer-order differential equations. However, proper initialization of fractional-order systems remains an ongoing problem [2-4]. This issue can date back as far as Riemann's complementary function theory, in which many mathematicians were made confused, including Liouville, Peacoch, Cayley, and Riemann himself [3].

In recent years, the initialization function approach and the infinite state approach provide two effective ways to impose physically coherent initial conditions to fractional systems. In the initialization function theory [3, 5], initialization functions are required to initialize fractional differential equations. The initialization function is a time-varying function and a generalization of the constant of integration required for the order-one integral. The Riemann-Liouville fractional derivatives and the Caputo derivatives can be initialized using the initialization function approach . In the infinite state theory, the Riemann-Liouville fractional integral with order between 0 and 1 is viewed a linear system which is characterized by the impulse function and excited by the integrand function. The linear system is termed as the fractional integrator, which is also an exactly equivalent frequency distributed differential system. As a result, initial conditions of fractional systems can be represented by the distributed initial conditions of the fractional integrator [6-8]. Equivalence and compatibility of the above initialization theories have been proved in [9, 10].

The Riemann-Liouville derivative and the Caputo derivative are the two most commonly used derivatives in the real-world modeling of factional systems. However, which derivative to choose is a trial and error process. In many applications, the two derivatives are expected to be equivalent as long as the initial conditions are properly taken into account [11, 12]. In [3] and [13], Laplace transforms of the Riemann-Liouville derivative and the Caputo derivative have been derived and shown to be equal in a special case where the history function is a constant and the fractional order lies between 0 and 1. In this paper, we present the equivalence and compatibility of the initialized Riemann-Liouville derivatives and the initialized Caputo derivatives with arbitrary fractional orders. By synthesizing the initialization function approach and the infinite state approach, the diffusive representations of the initialized Riemann-Liouville derivatives and the initialized Caputo derivatives are derived. Laplace transforms of the two initialized derivatives are shown to be equal. Consequently, the Riemann-Liouville derivatives and the Caputo derivatives are proved to equivalent when initial conditions are properly imposed.

The rest of this paper is organized as follows. Section 2 revisits the diffusive representation for the initialized Riemann-Liouville fractional integral. Section 3 presents the equivalence of the initialized Riemann-Liouville derivatives and the initialized Caputo derivatives with order between 0 and 1. Section 4 shows the equivalence of the two initialized derivatives with order between 1 and 2. Section 5 proves the equivalence of the two initialized derivatives with arbitrary orders. Finally, the paper is concluded in Section 6.

\section{Preliminaries}
\subsection{Diffusive model of the fractional integrator [7]}
The Riemann-Liouville fractional integral of a function $f\left( t \right)$with order $0<\alpha <1$ is defined as
\begin{equation}
{}_{{{t}_{0}}}I_{t}^{\alpha }f\left( t \right)=\frac{1}{\Gamma \left( \alpha  \right)}\int_{{{t}_{0}}}^{t}{{{\left( t-\tau  \right)}^{\alpha -1}}}f\left( \tau  \right)d\tau
\end{equation}
Eq.(1)can also be viewed as a convolution between $f\left( t \right)$ and the impulse response
${{h}_{\alpha }}\left( t \right)=\frac{{{t}^{\alpha -1}}}{\text{ }\!\!\Gamma\!\!\text{ }\left( \alpha  \right)}$, i.e.,
$${}_{{{t}_{0}}}I_{t}^{\alpha }f(t)={{h}_{\alpha }}\left( t \right)*f\left( t \right).$$
From this viewpoint, the fractional integral can be obtained as the output of a linear system, which is characterized by the impulse response ${{h}_{\alpha }}\left( t \right)$ and excited by $f\left( t \right)$:
$$x(t)={{h}_{\alpha }}\left( t \right)*f\left( t \right).$$
This linear system is called the fractional integrator. The frequency distributed state $z(\omega ,t)$of the fractional integrator verifies the following ordinary differential equation:
\begin{equation}
\frac{\partial z(\omega ,t)}{\partial t}=-\omega z(\omega ,t)+f(t)
\end{equation}
and the output of the fractional integrator is
\begin{equation}
x(t)=\int_{0}^{+\infty }{{{\mu }_{\alpha }}\left( \omega  \right)z(\omega ,t)}d\omega,
\end{equation}
where the elementary frequency $\omega $ is ranging from 0 to $\infty $, and
\[{{\mu }_{\alpha }}(\omega )=\frac{sin\alpha \pi }{\pi }{{\omega }^{-\alpha }}.\]
The relations (2) and (3) are called the frequency distributed model or diffusive model of the fractional integrator.

\subsection{The initialized fractional integrals [14, 15]}
In the initialization theory for fractional calculus, the integrand of Riemann-Liouville fractional integral is represented as a composite function
\[v(t)=\left\{ \begin{aligned}
  & 0,\text{  }t<-a \\
 & {{f}_{in}}(t),\text{  }-a\le t\le 0 \\
 & f(t),\text{  }t>0 \\
\end{aligned} \right.\]
where $t=-a$ is the starting time of integral, $t=0$ is the initial time, ${{f}_{i}}(t)$ is the function describing the behavior during the initialization period $\left[ -a,0 \right]$, $f(t)$ is the function of primary interest after the initial time $t=0$.

The initialized Riemann-Liouville fractional integral of order $\alpha $ is defined as
$${}_{0}D_{t}^{-\alpha }f\left( t \right)=\frac{1}{\Gamma \left( \alpha  \right)}\int_{0}^{t}{{{\left( t-\tau  \right)}^{\alpha -1}}}f\left( \tau  \right)d\tau +\psi \left( t \right),\text{ }t>0$$
or
$${}_{0}D_{t}^{-\alpha }f\left( t \right)={}_{0}d_{t}^{-\alpha }f\left( t \right)+\psi \left( t \right),\text{ }t>0.$$
where ${}_{0}D_{t}^{-\alpha }f\left( t \right)$ is called the initialized fractional integral, ${}_{0}d_{t}^{-\alpha }f\left( t \right)$ is called the uninitialized fractional integral, $\psi \left( t \right)$ is called the initialization function describing the hereditary effect of the past.
\subsection{Diffusive representation of the initialized Riemann-Liouville fractional integrals [10]}
\begin{lemma}
\newtheorem{lem}{Lemma} The uninitialized Riemann-Liouville fractional integral ${}_{0}d_{t}^{-\alpha }f(t)$ with order $0<\alpha <1$ is equivalent to the diffusive model with zero initial condition:
\[\left\{ \begin{aligned}
  & \frac{\partial z(\omega ,t)}{\partial t}=-\omega z(\omega ,t)+f(t) \\
 & z(\omega ,{{t}_{0}})=0 \\
\end{aligned} \right.\]
and
\[{}_{0}d_{t}^{-\alpha }f(t)=\int_{0}^{+\infty }{{{\mu }_{\alpha }}\left( \omega  \right)z(\omega ,t)}d\omega. \]
\end{lemma}
\begin{lemma} The initialized Riemann-Liouville fractional integral ${}_{0}D_{t}^{-\alpha }f(t)$ with order $0<\alpha <1$ is equivalent to the diffusive model with distributed initial condition:
	\[\left\{ \begin{aligned}
  & \frac{\partial z(\omega ,t)}{\partial t}=-\omega z(\omega ,t)+f(t),\text{  } \\
 & z(\omega ,0)=\int_{-a}^{0}{{{e}^{-\omega \tau }}{{f}_{in}}(\tau )d\tau } \\
\end{aligned} \right.\]	
and
\[{}_{0}D_{t}^{-\alpha }f(t)=\int_{0}^{+\infty }{{{\mu }_{\alpha }}(\omega )z(\omega ,t)d\omega }.\]
\end{lemma}
\section{Equivalence of the two Derivatives with order between 0 and 1}
\subsection{Diffusive model for the initialized Riemann-Liouville derivatives with order between 0 and 1}
The initialized Riemann-Liouville fractional derivative ${}_{0}D_{t}^{\alpha }f(t)$ with order $0<\alpha <1 $
is defined in [14] as
$${}_{0}^{RL}D_{t}^{\alpha }f\left( t \right)=\frac{d}{dt}{}_{0}D_{t}^{\alpha-1 }f\left( t \right).$$
By virtue of Lemma 2, the diffusive representation of the initialized Riemann-Liouville fractional integral ${}_{0}D_{t}^{\alpha-1 }f(t)$ is
$${}_{0}D_{t}^{\alpha-1}f\left( t \right)=\int_{0}^{\infty }{{{\mu }_{1\text{-}\alpha }}(\omega ){{z}_{RL}}\left( \omega ,t \right)d\omega }$$
where ${{z}_{RL}}\left( \omega ,t \right)$ satisfies
\begin{equation}
\begin{cases}
\frac{\partial {{z}_{RL}}(\omega ,t)}{\partial t}=-\omega {{z}_{RL}}(\omega ,t)+f(t) \\
{{z}_{RL}}(\omega ,0)=\int_{-a}^{0}{{{e}^{-\omega \tau }}{{f}_{in}}(\tau )d\tau } \\
\end{cases}
\end{equation}
As a result, the diffusive model for the initialized Riemann-Liouville derivatives is
\begin{equation}
{}_{0}^{RL}D_{t}^{\alpha }f\left( t \right)=\frac{d}{dt}\int_{0}^{\infty }{{{\mu }_{1\text{-}\alpha }}(\omega ){{z}_{RL}}\left( \omega ,t \right)d\omega }
\end{equation}
Taking Laplace transform of the first equation of Eq.(4), we have
\begin{equation}
{{Z}_{RL}}\left( \omega ,t \right)=\frac{{{z}_{RL}}\left( \omega ,0 \right)+F\left( s \right)}{s+\omega }
\end{equation}
Taking Laplace transform of the initialized Riemann-Liouville derivative from Eq.(5), we have
\begin{equation}
\begin{split}
L\left\{ {}_{0}^{RL}D_{t}^{\alpha }f\left( t \right) \right\}&=L\left\{ \frac{d}{dt}\int_{0}^{\infty }{{{\mu }_{1\text{-}\alpha }}(\omega ){{z}_{RL}}\left( \omega ,t \right)d\omega } \right\} \\
 & =s\int_{0}^{\infty }{{{\mu }_{1\text{-}\alpha }}(\omega ){{Z}_{RL}}\left( \omega ,t \right)d\omega }-\int_{0}^{\infty }{{{\mu }_{1\text{-}\alpha }}(\omega ){{z}_{RL}}\left( \omega ,0 \right)d\omega } \\
\end{split}
\end{equation}
Substituting Eq. (6) into Eq.(7), we obtain
\begin{equation}
\begin{split}
L\left\{ {}_{0}^{RL}D_{t}^{\alpha }f\left( t \right) \right\}&=sF\left( s \right)\int_{0}^{\infty }{\frac{{{\mu }_{1\text{-}\alpha }}(\omega )}{s+\omega }d\omega }+s\int_{0}^{\infty }{\frac{{{\mu }_{1\text{-}\alpha }}(\omega ){{z}_{RL}}\left( \omega ,0 \right)}{s+\omega }d\omega } \\
 & \quad -\int_{0}^{\infty }{{{\mu }_{1\text{-}\alpha }}(\omega ){{z}_{RL}}\left( \omega ,0 \right)d\omega } \\
\end{split}
\end{equation}
Because
$$\int_{0}^{\infty }{\frac{{{\mu }_{\alpha }}(\omega )}{s+\omega }d\omega }={{s}^{-\alpha }},$$
Eq.(8) becomes
\begin{equation}
\begin{split}
L\left\{ {}_{0}^{RL}D_{t}^{\alpha }f\left( t \right) \right\}&={{s}^{\alpha }}F\left( s \right)+s\int_{0}^{\infty }{\frac{{{\mu }_{1\text{-}\alpha }}(\omega ){{z}_{RL}}\left( \omega ,0 \right)}{s+\omega }d\omega }\\
&\quad -\int_{0}^{\infty }{{{\mu }_{1\text{-}\alpha }}(\omega ){{z}_{RL}}\left( \omega ,0 \right)d\omega } \\
 & ={{s}^{\alpha }}F\left( s \right)+\int_{0}^{\infty }{\left( \frac{s}{s+\omega }-1 \right){{\mu }_{1\text{-}\alpha }}(\omega ){{z}_{RL}}\left( \omega ,0 \right)d\omega } \\
 & ={{s}^{\alpha }}F\left( s \right)-\int_{0}^{\infty }{\frac{\omega {{\mu }_{1\text{-}\alpha }}(\omega ){{z}_{RL}}\left( \omega ,0 \right)}{s+\omega }d\omega } \\
\end{split}
\end{equation}
\subsection{Diffusive model for initialized Caputo derivatives with order between 0 and 1}
The initialized Caputo fractional derivative ${}_{0}D_{t}^{\alpha }f(t)$ with order $0<\alpha <1$ is defined in [14] as
                 $${}_{0}^{C}D_{t}^{\alpha }f\left( t \right)={}_{0}D_{t}^{\alpha-1}{f}'\left( t \right).$$
By virtue of Lemma 2, the diffusive representation of the initialized Caputo fractional derivative is
\begin{equation}
{}_{0}^{C}D_{t}^{\alpha }f\left( t \right)=\int_{0}^{\infty }{{{\mu }_{1\text{-}\alpha }}(\omega ){{z}_{C}}\left( \omega ,t \right)d\omega }
\end{equation}
where ${{z}_{C}}\left( \omega ,t \right)$ satisfies

\begin{equation}
\begin{cases}
\frac{\partial {{z}_{C}}(\omega ,t)}{\partial t}=-\omega {{z}_{C}}(\omega ,t)+f(t) \\
{{z}_{C}}(\omega ,0)=\int_{-a}^{0}{{{e}^{-\omega \tau }}{{f}_{in}}^{\prime }(\tau )d\tau } \\
\end{cases}
\end{equation}
Taking Laplace transform of Eq.(11), we have
\begin{equation}
{{Z}_{C}}\left( \omega ,t \right)=\frac{sF\left( s \right)-f\left( 0 \right)+{{z}_{C}}\left( \omega ,0 \right)}{s+\omega }
\end{equation}	
Taking Laplace transform of Eq.(10), we have
\begin{equation}
\begin{split}
L\left\{ {}_{0}^{C}D_{t}^{\alpha }f\left( t \right) \right\}&=L\left\{ \int_{0}^{\infty }{{{\mu }_{1\text{-}\alpha }}(\omega ){{z}_{C}}\left( \omega ,t \right)d\omega } \right\} \\
 & =\int_{0}^{\infty }{{{\mu }_{1\text{-}\alpha }}(\omega ){{Z}_{C}}\left( \omega ,t \right)d\omega } \\
\end{split}
\end{equation}
Substituting Eq.(12) into Eq.(13), we obtain
\begin{equation}
\begin{split}
L\left\{ {}_{0}^{C}D_{t}^{\alpha }f\left( t \right) \right\}&=sF\left( s \right)\int_{0}^{\infty }{\frac{{{\mu }_{1\text{-}\alpha }}(\omega )}{s+\omega }d\omega }+\int_{0}^{\infty }{\frac{{{\mu }_{1\text{-}\alpha }}(\omega ){{z}_{C}}\left( \omega ,0 \right)}{s+\omega }d\omega } \\
 & \quad -f\left( 0 \right)\int_{0}^{\infty }{\frac{{{\mu }_{1\text{-}\alpha }}(\omega )}{s+\omega }d\omega } \\
 & ={{s}^{\alpha }}F\left( s \right)+\int_{0}^{\infty }{\frac{{{\mu }_{1\text{-}\alpha }}(\omega ){{z}_{C}}\left( \omega ,0 \right)}{s+\omega }d\omega }-{{s}^{\alpha -1}}f\left( 0 \right) \\
\end{split}
\end{equation}
In terms of Eq.(11) and integrating by parts, we have
\begin{equation*}
{{z}_{C}}(\omega ,0)=\int_{-a}^{0}{{{e}^{-\omega \tau }}{{f}_{in}}^{\prime }(\tau )d\tau }=\left. {{e}^{-\omega \tau }}{{f}_{in}}(\tau ) \right|_{\tau =-a}^{\tau =0}-\omega \int_{-a}^{0}{{{e}^{-\omega \tau }}{{f}_{in}}(\tau )d\tau }
\end{equation*}
In terms of Eq.(4), we have

\begin{equation}
{{z}_{C}}(\omega ,0)=f\left( 0 \right)-\omega {{z}_{Rl}}(\omega ,0)
\end{equation}
Then
\begin{equation}
\begin{split}
\int_{0}^{\infty }{\frac{{{\mu }_{1\text{-}\alpha }}(\omega ){{z}_{C}}\left( \omega ,0 \right)}{s+\omega }d\omega } & =f\left( 0 \right)\int_{0}^{\infty }{\frac{{{\mu }_{1-\alpha }}(\omega )}{s+\omega }d\omega }-\int_{0}^{\infty }{\frac{\omega {{\mu }_{1\text{-}\alpha }}(\omega ){{z}_{RL}}\left( \omega ,0 \right)}{s+\omega }d\omega } \\
 & ={{s}^{\alpha -1}}f\left( 0 \right)-\int_{0}^{\infty }{\frac{{{\mu }_{1\text{-}\alpha }}(\omega ){{z}_{RL}}\left( \omega ,0 \right)}{s+\omega }d\omega } \\
\end{split}
\end{equation}
Substituting Eq.(16) into Eq.(14), we get
\begin{equation}
L\left\{ {}_{0}^{C}D_{t}^{\alpha }f\left( t \right) \right\}={{s}^{\alpha }}F\left( s \right)-\int_{0}^{\infty }{\frac{\omega {{\mu }_{1\text{-}\alpha }}(\omega ){{z}_{RL}}\left( \omega ,0 \right)}{s+\omega }d\omega }
\end{equation}
From Eq.(9) and Eq.(17),  it is clear that
\begin{equation}
L\left\{ {}_{0}^{RL}D_{t}^{\alpha }f\left( t \right) \right\}=L\left\{ {}_{0}^{C}D_{t}^{\alpha }f\left( t \right) \right\}.
\end{equation}
Thus
\begin{equation}
{}_{0}^{RL}D_{t}^{\alpha }f\left( t \right)\text{=}{}_{0}^{C}D_{t}^{\alpha }f\left( t \right).
\end{equation}
Eq.(19) shows the equivalence of the initialized Riemann-Liouville derivatives and Caputo derivatives with order between 0 and 1.

\section{Equivalence of the two Derivatives with order between 1 and 2}
\subsection{Diffusive model for the initialized Riemann-Liouville derivatives with order between 1 and 2}
The initialized Riemann-Liouville fractional derivative ${}_{0}D_{t}^{\alpha }f(t)$ with order $1<\alpha <2$ is defined as
$${}_{0}^{RL}D_{t}^{\alpha }f\left( t \right)=\frac{{{d}^{2}}}{d{{t}^{2}}}{}_{0}D_{t}^{\alpha-2}f\left( t \right).$$
By virtue of Lemma 2, the diffusive representation of the initialized Riemann-Liouville fractional integral ${}_{0}D_{t}^{\alpha-2}f(t)$ is
$${}_{0}D_{t}^{\alpha-2}f\left( t \right)=\int_{0}^{\infty }{{{\mu }_{2-\alpha }}(\omega ){{z}_{RL}}\left( \omega ,t \right)d\omega }$$
where ${{z}_{RL}}\left( \omega ,t \right)$ satisfies Eq.(4).
As a result, the diffusive model for the initialized Riemann-Liouville derivatives is
\begin{equation}
{}_{0}^{RL}D_{t}^{\alpha }f\left( t \right)=\frac{{{d}^{2}}}{d{{t}^{2}}}\int_{0}^{\infty }{{{\mu }_{1\text{-}\alpha }}(\omega ){{z}_{RL}}\left( \omega ,t \right)d\omega }
\end{equation}
Taking Laplace transform of Eq.(20), we have
\begin{equation}
\begin{split}
L\left\{ {}_{0}^{RL}D_{t}^{\alpha }f\left( t \right) \right\}&=L\left\{ \frac{{{d}^{2}}}{d{{t}^{2}}}\int_{0}^{\infty }{{{\mu }_{1\text{-}\alpha }}(\omega ){{z}_{RL}}\left( \omega ,t \right)d\omega } \right\} \\
 & ={{s}^{2}}\int_{0}^{\infty }{{{\mu }_{\text{2-}\alpha }}(\omega ){{Z}_{RL}}\left( \omega ,t \right)d\omega }-s\int_{0}^{\infty }{{{\mu }_{\text{2-}\alpha }}(\omega ){{z}_{RL}}\left( \omega ,0 \right)d\omega } \\
 & \quad-\int_{0}^{\infty }{{{\mu }_{\text{2-}\alpha }}(\omega ){{\left. \frac{\partial {{z}_{RL}}(\omega ,t)}{\partial t} \right|}_{t=0}}d\omega } \\
\end{split}
\end{equation}
Substituting Eq.(6) into Eq.(21), we have
\begin{equation}
\begin{split}
L\left\{ {}_{0}^{RL}D_{t}^{\alpha }f\left( t \right) \right\}&=sF\left( s \right)\int_{0}^{\infty }{\frac{{{\mu }_{\text{2-}\alpha }}(\omega )}{s+\omega }d\omega }+{{s}^{2}}\int_{0}^{\infty }{\frac{{{\mu }_{\text{2-}\alpha }}(\omega ){{z}_{RL}}\left( \omega ,0 \right)}{s+\omega }d\omega } \\
 & \quad-s\int_{0}^{\infty }{{{\mu }_{\text{2-}\alpha }}(\omega ){{z}_{RL}}\left( \omega ,0 \right)d\omega -\int_{0}^{\infty }{{{\mu }_{2-\alpha }}(\omega ){{\left. \frac{\partial {{z}_{RL}}(\omega ,t)}{\partial t} \right|}_{t=0}}d\omega }} \\
\end{split}
\end{equation}
Because
$$\int_{0}^{\infty }{\frac{{{\mu }_{2-\alpha }}(\omega )}{s+\omega }d\omega }={{s}^{\alpha -2}},$$
By substituting the first equation of Eq.(4) into Eq.(22), we obtain
\begin{equation}
\begin{split}
L\left\{ {}_{0}^{RL}D_{t}^{\alpha }f\left( t \right) \right\}&={{s}^{\alpha }}F\left( s \right)+\int_{0}^{\infty }{\left( \frac{{{s}^{2}}}{s+\omega }-s+\omega  \right){{\mu }_{\text{2-}\alpha }}(\omega ){{z}_{RL}}\left( \omega ,0 \right)d\omega } \\
 & \quad-f\left( 0 \right)\int_{0}^{\infty }{{{\mu }_{2-\alpha }}(\omega )d\omega } \\
 & ={{s}^{\alpha }}F\left( s \right)+\int_{0}^{\infty }{\frac{{{\omega }^{2}}{{\mu }_{\text{2-}\alpha }}(\omega ){{z}_{RL}}\left( \omega ,0 \right)}{s+\omega }d\omega }\\
 & \quad-f\left( 0 \right)\int_{0}^{\infty }{{{\mu }_{\text{2-}\alpha }}(\omega )d\omega } \\
\end{split}
\end{equation}

\subsection{Diffusive model for the initialized Caputo derivatives with order between 1 and 2}
The initialized Caputo fractional derivative ${}_{0}D_{t}^{\alpha }f(t)$ with order $1<\alpha <2$ is defined as
                 $${}_{0}^{C}D_{t}^{\alpha }f\left( t \right)={}_{0}D_{t}^{\alpha-2 }{f}''\left( t \right)$$
By virtue of Lemma 2, the diffusive representation of the initialized Caputo fractional derivative is
\begin{equation}
{}_{0}^{C}D_{t}^{\alpha }f\left( t \right)=\int_{0}^{\infty }{{{\mu }_{\text{2-}\alpha }}(\omega ){{z}_{C}}\left( \omega ,t \right)d\omega }
\end{equation}
where ${{z}_{C}}\left( \omega ,t \right)$ satisfies
\begin{equation}
\begin{cases}
\frac{\partial {{z}_{C}}(\omega ,t)}{\partial t}=-\omega {{z}_{C}}(\omega ,t)+{f}''\left( t \right) \\
{{z}_{C}}(\omega ,0)=\int_{-a}^{0}{{{e}^{-\omega \tau }}{{f}_{in}}^{\prime \prime }(\tau )d\tau } \\
\end{cases}
\end{equation}
Taking Laplace transform of the first equation in Eq.(25), we have
\begin{equation}
{{Z}_{C}}\left( \omega ,t \right)=\frac{{{s}^{2}}F\left( s \right)-sf\left( 0 \right)-{f}'\left( 0 \right)+{{z}_{C}}\left( \omega ,0 \right)}{s+\omega }.
\end{equation}
In terms of the initial conditions of Eq.(25),

\begin{equation}
\begin{split}
{{z}_{C}}(\omega ,0)&=\int_{-a}^{0}{{{e}^{-\omega \tau }}{{f}_{in}}^{\prime \prime }(\tau )d\tau } \\
 & =\left. {{e}^{-\omega \tau }}{{f}_{in}}^{\prime }(\tau ) \right|_{\tau =-a}^{\tau =0}-\omega \int_{-a}^{0}{{{e}^{-\omega \tau }}{{f}_{in}}^{\prime }(\tau )d\tau } \\
 & ={f}'\left( 0 \right)-\omega \int_{-a}^{0}{{{e}^{-\omega \tau }}{{f}_{in}}^{\prime }(\tau )d\tau } \\
 & ={f}'\left( 0 \right)-\omega f\left( 0 \right)+{{\omega }^{2}}{{z}_{R}}(\omega ,0) \\
\end{split}
\end{equation}
Substituting Eq.(27) into Eq.(26), we have
\begin{equation}
{{Z}_{C}}\left( \omega ,t \right)=\frac{{{s}^{2}}F\left( s \right)-\left( s+\omega  \right)f\left( 0 \right)+{{\omega }^{2}}{{z}_{RL}}\left( \omega ,0 \right)}{s+\omega }.
\end{equation}
Taking Laplace transform of Eq.(24)and substituting Eq.(28) into it, we obtain

\begin{equation}
\begin{split}
L\left\{ {}_{0}^{C}D_{t}^{\alpha }f\left( t \right) \right\}&=L\left\{ \int_{0}^{\infty }{{{\mu }_{\text{2-}\alpha }}(\omega ){{z}_{C}}\left( \omega ,t \right)d\omega } \right\} \\
 & =\int_{0}^{\infty }{{{\mu }_{\text{2-}\alpha }}(\omega ){{Z}_{C}}\left( \omega ,t \right)d\omega } \\
 & ={{s}^{2}}F\left( s \right)\int_{0}^{\infty }{\frac{{{\mu }_{2-\alpha }}(\omega )}{s+\omega }d\omega }+\int_{0}^{\infty }{\frac{{{\omega }^{2}}{{\mu }_{\text{2-}\alpha }}(\omega ){{z}_{RL}}\left( \omega ,0 \right)}{s+\omega }d\omega }\\
 &\quad -f\left( 0 \right)\int_{0}^{\infty }{{{\mu }_{\text{2-}\alpha }}(\omega )d\omega } \\
 & ={{s}^{\alpha }}F\left( s \right)+\int_{0}^{\infty }{\frac{{{\omega }^{2}}{{\mu }_{\text{2-}\alpha }}(\omega ){{z}_{RL}}\left( \omega ,0 \right)}{s+\omega }d\omega }\\
 &\quad -f\left( 0 \right)\int_{0}^{\infty }{{{\mu }_{\text{2-}\alpha }}(\omega )d\omega } \\
\end{split}
\end{equation}
From Eq.(23) and Eq.(29), it is clear that
\begin{equation}
L\left\{ {}_{0}^{RL}D_{t}^{\alpha }f\left( t \right) \right\}=L\left\{ {}_{0}^{C}D_{t}^{\alpha }f\left( t \right) \right\}.
\end{equation}
Thus,
\begin{equation}
{}_{0}^{RL}D_{t}^{\alpha }f\left( t \right)\text{=}{}_{0}^{C}D_{t}^{\alpha }f\left( t \right).
\end{equation}
Eq.(31) shows the equivalence of the initialized Riemann-Liouville derivatives and Caputo derivatives with order between 1 and 2.
\section{Equivalence of the two Derivatives with arbitrary orders}
\subsection{Diffusive model for initialized Riemann-Liouville derivatives with order between $n$ and $n-1$}
The initialized Riemann-Liouville fractional derivative ${}_{0}D_{t}^{\alpha }f(t)$ with order $n<\alpha <n-1$
is defined as
$${}_{0}^{RL}D_{t}^{\alpha }f\left( t \right)=\frac{{{d}^{n}}}{d{{t}^{n}}}{}_{0}D_{t}^{\alpha-n}f\left( t \right).$$
By virtue of Lemma 2, the diffusive representation of the initialized Riemann-Liouville fractional integral ${}_{0}D_{t}^{\alpha-n }f\left( t \right)$ is
$${}_{0}D_{t}^{\alpha-n}f\left( t \right)=\int_{0}^{\infty }{{{\mu }_{n-\alpha }}(\omega ){{z}_{RL}}\left( \omega ,t \right)d\omega }$$
where ${{z}_{RL}}\left( \omega ,t \right)$ satisfies Eq.(4).

As a result, the diffusive model for the initialized Riemann-Liouville derivative is
\begin{equation}
{}_{0}^{RL}D_{t}^{\alpha }f\left( t \right)=\frac{{{d}^{n}}}{d{{t}^{n}}}\int_{0}^{\infty }{{{\mu }_{\text{n-}\alpha }}(\omega ){{z}_{RL}}\left( \omega ,t \right)d\omega }
\end{equation}
Taking Laplace transform of Eq.(32), we have

\begin{equation}
\begin{split}
L\left\{ {}_{0}^{RL}D_{t}^{\alpha }f\left( t \right) \right\}&=L\left\{ \frac{{{d}^{n}}}{d{{t}^{n}}}\int_{0}^{\infty }{{{\mu }_{n-\alpha }}(\omega ){{z}_{RL}}\left( \omega ,t \right)d\omega } \right\} \\
 & ={{s}^{n}}\int_{0}^{\infty }{{{\mu }_{n-\alpha }}(\omega ){{Z}_{RL}}\left( \omega ,t \right)d\omega }-{{s}^{n-1}}\int_{0}^{\infty }{{{\mu }_{n-\alpha }}(\omega ){{z}_{RL}}\left( \omega ,0 \right)d\omega } \\
 & \quad-{{s}^{n-2}}\int_{0}^{\infty }{{{\mu }_{n-\alpha }}(\omega ){{{{z}'}}_{RL}}(\omega ,0)d\omega } \\
 & \quad-{{s}^{n-3}}\int_{0}^{\infty }{{{\mu }_{n-\alpha }}(\omega ){{{{z}''}}_{RL}}(\omega ,0)d\omega } \\
 & \quad-\cdots -\int_{0}^{\infty }{{{\mu }_{n-\alpha }}(\omega ){{z}^{\left( n-1 \right)}}_{RL}(\omega ,0)d\omega } \\
\end{split}
\end{equation}
Substituting Eq.(6) into Eq.(33), and applying $\int_{0}^{\infty }{\frac{{{\mu }_{2-\alpha }}(\omega )}{s+\omega }d\omega }={{s}^{\alpha -2}}$, we have

\begin{equation}
\begin{split}
L\left\{ {}_{0}^{RL}D_{t}^{\alpha }f\left( t \right) \right\}&={{s}^{\alpha }}F\left( s \right)+{{s}^{n}}\int_{0}^{\infty }{\frac{{{\mu }_{n-\alpha }}(\omega ){{z}_{RL}}\left( \omega ,0 \right)}{s+\omega }d\omega } \\
 & \quad -{{s}^{n-1}}\int_{0}^{\infty }{{{\mu }_{n-\alpha }}(\omega ){{z}_{RL}}\left( \omega ,0 \right)d\omega } \\
 & \quad -{{s}^{n-2}}\int_{0}^{\infty }{{{\mu }_{n-\alpha }}(\omega ){{{{z}'}}_{RL}}\left( \omega ,0 \right)d\omega } \\
 & \quad -{{s}^{n-3}}\int_{0}^{\infty }{{{\mu }_{n-\alpha }}(\omega ){{{{z}''}}_{RL}}\left( \omega ,0 \right)d\omega } \\
 & \quad -\cdots -\int_{0}^{\infty }{{{\mu }_{n-\alpha }}(\omega )z_{RL}^{^{\left( n-1 \right)}}\left( \omega ,0 \right)d\omega } \\
\end{split}
\end{equation}

From Eq.(4), we can calculate high-order derivatives of ${{z}_{RL}}\left( \omega ,t \right)$:
$$\frac{{{\partial }^{2}}{{z}_{RL}}\left( \omega ,t \right)}{\partial {{t}^{2}}}={{\omega }^{2}}{{z}_{RL}}\left( \omega ,t \right)-\omega f\left( t \right)+{f}''\left( t \right)$$
$$\frac{{{\partial }^{3}}{{z}_{RL}}\left( \omega ,t \right)}{\partial {{t}^{3}}}=-{{\omega }^{3}}{{z}_{RL}}\left( \omega ,t \right)+{{\omega }^{2}}f\left( t \right)-\omega {f}'\left( t \right)+{f}''\left( t \right)$$
$$\cdots \cdots$$
\begin{equation*}\begin{split}
\frac{{{\partial }^{n-1}}{{z}_{RL}}\left( \omega ,t \right)}{\partial {{t}^{n-1}}}&={{\left( -\omega  \right)}^{n-1}}{{z}_{RL}}\left( \omega ,t \right)+{{\left( -\omega  \right)}^{n-2}}f\left( t \right)+{{\left( -\omega  \right)}^{n-3}}{f}'\left( t \right)\\
& \quad +\cdots +{{f}^{\left( n-2 \right)}}\left( t \right)
\end{split}\end{equation*}
Substituting the above derivatives into Eq.(34), we obtain

\begin{equation}
\begin{split}
L\left\{ {}_{0}^{RL}D_{t}^{\alpha }f\left( t \right) \right\}&={{s}^{\alpha }}F\left( s \right)+{{s}^{n}}\int_{0}^{\infty }{\frac{{{\mu }_{n-\alpha }}(\omega ){{z}_{RL}}\left( \omega ,0 \right)}{s+\omega }d\omega } \\
 &  \quad-{{s}^{n-1}}\int_{0}^{\infty }{{{\mu }_{n-\alpha }}(\omega ){{z}_{RL}}\left( \omega ,0 \right)d\omega } \\
 &  \quad-{{s}^{n-2}}\int_{0}^{\infty }{{{\mu }_{n-\alpha }}(\omega )\left[ -\omega {{z}_{RL}}(\omega ,0)+f(0) \right]d\omega } \\
 &  \quad-{{s}^{n-3}}\int_{0}^{\infty }{{{\mu }_{n-\alpha }}(\omega )\left[ {{\omega }^{2}}{{z}_{RL}}\left( \omega ,0 \right)-\omega f\left( 0 \right)+{f}''\left( 0 \right) \right]d\omega } \\
 &  \quad-\cdots \cdots  \\
 &  \quad-\int_{0}^{\infty }{{{\mu }_{n-\alpha }}(\omega )\left[ {{\left( -\omega  \right)}^{n-1}}{{z}_{RL}}\left( \omega ,0 \right)+{{\left( -\omega  \right)}^{n-2}}f\left( 0 \right) \right.} \\
 &  \quad\left. +{{\left( -\omega  \right)}^{n-3}}{f}'\left( 0 \right)+\cdots +{{f}^{\left( n-2 \right)}}\left( 0 \right) \right]d\omega  \\
 & ={{s}^{\alpha }}F\left( s \right)+\int_{0}^{\infty }{{{\mu }_{n-\alpha }}(\omega ){{z}_{RL}}\left( \omega ,0 \right)\left[ \frac{{{s}^{n}}}{s+\omega }-{{s}^{n-1}}+\omega {{s}^{n-2}} \right.} \\
 &  \quad\left. -{{\omega }^{2}}{{s}^{n-3}}+\cdots +{{\left( -\omega  \right)}^{n-2}}s+{{\left( -\omega  \right)}^{n-1}} \right]d\omega \\
 &  \quad +\int_{0}^{\infty }{{{\mu }_{n-\alpha }}(\omega )\left\{ {{\Delta }_{1}} \right\}d\omega } \\
 & ={{s}^{\alpha }}F\left( s \right)+\int_{0}^{\infty }{{{\mu }_{n-\alpha }}(\omega ){{z}_{RL}}\left( \omega ,0 \right)\frac{{{\left( -\omega  \right)}^{n}}}{s+\omega }d\omega } \\
 &  \quad-\int_{0}^{\infty }{{{\mu }_{n-\alpha }}(\omega )\left\{ {{\Delta }_{1}} \right\}d\omega } \\
\end{split}
\end{equation}
where
\begin{equation}
\begin{split}
{{\Delta }_{1}}&=\left[ {{s}^{n-2}}+\left( -\omega  \right){{s}^{n-3}}+{{\left( -\omega  \right)}^{2}}{{s}^{n-4}}+\cdots {{\left( -\omega  \right)}^{n-2}} \right]f\left( 0 \right) \\
 & \quad+\left[ {{s}^{n-3}}+\left( -\omega  \right){{s}^{n-4}}+\cdots {{\left( -\omega  \right)}^{n-3}} \right]{f}'\left( 0 \right)+\cdots +{{f}^{\left( n-2 \right)}}\left( 0 \right) \\
\end{split}
\end{equation}
\subsection{Diffusive model for initialized Caputo derivatives with order between $n$ and $n-1$}
The initialized Caputo fractional derivative ${}_{0}D_{t}^{\alpha }f(t)$ with order $n-1<\alpha <n$ is defined as
                 $${}_{0}^{C}D_{t}^{\alpha }f\left( t \right)={}_{0}D_{t}^{\alpha-n}{{f}^{\left( n \right)}}\left( t \right).$$

By virtue of Lemma 2, the diffusive representation of the initialized Caputo fractional derivative is
\begin{equation}
{}_{0}^{C}D_{t}^{\alpha }f\left( t \right)=\int_{0}^{\infty }{{{\mu }_{n-\alpha }}(\omega ){{z}_{C}}\left( \omega ,t \right)d\omega }
\end{equation}
where ${{z}_{C}}\left( \omega ,t \right)$ satisfies
\begin{equation}
\begin{cases}
\frac{\partial {{z}_{C}}(\omega ,t)}{\partial t}=-\omega {{z}_{C}}(\omega ,t)+{{f}^{\left( n \right)}}\left( t \right) \\
{{z}_{C}}(\omega ,0)=\int_{-a}^{0}{{{e}^{-\omega \tau }}f_{in}^{\left( n \right)}(\tau )d\tau } \\
\end{cases}
\end{equation}

Taking Laplace transform of Eq.(38), we have
\begin{equation}
\begin{split}
{{Z}_{C}}\left( \omega ,s \right)&=\frac{1}{s+\omega }\left[ {{z}_{C}}\left( \omega ,0 \right)+{{s}^{n}}F\left( s \right)-{{s}^{n-1}}f\left( 0 \right)-{{s}^{n-1}}{f}'\left( 0 \right)-\cdots  \right. \\
 &  \quad \left. -s{{f}^{\left( n-2 \right)}}\left( 0 \right)-{{f}^{\left( n-1 \right)}}\left( 0 \right) \right] \\
\end{split}
\end{equation}
In terms of the initial condition of Eq.(38) and integrating by parts, we have
\begin{equation}
\begin{split}
{{z}_{C}}(\omega ,0)&={{f}^{\left( n-1 \right)}}\left( 0 \right)+\left( -\omega  \right){{f}^{\left( n-2 \right)}}\left( 0 \right)+{{\left( -\omega  \right)}^{2}}{{f}^{\left( n-3 \right)}}\left( 0 \right)+{{\left( -\omega  \right)}^{3}}{{f}^{\left( n-4 \right)}}\left( 0 \right) \\
 & \quad +\cdots +{{\left( -\omega  \right)}^{n-2}}{f}'\left( 0 \right)+{{\left( -\omega  \right)}^{n-1}}f\left( 0 \right)+{{\left( -\omega  \right)}^{n}}{{z}_{RL}}(\omega ,0) \\
\end{split}
\end{equation}
Substituting Eq.(40) into Eq.(39), we have
\begin{equation}
\begin{split}
{{Z}_{C}}\left( \omega ,s \right)&=\frac{1}{s+\omega }\left\{ {{s}^{n}}F\left( s \right)+\left[ {{\left( -\omega  \right)}^{n-1}}-{{s}^{n-1}} \right]f\left( 0 \right) \right. \\
 & \quad+\left[ {{\left( -\omega  \right)}^{n-2}}-{{s}^{n-2}} \right]{f}'\left( 0 \right)+\cdots  \\
 & \quad\left. +\left[ \left( -\omega  \right)-s \right]{{f}^{\left( n-2 \right)}}\left( 0 \right)+{{\left( -\omega  \right)}^{n}}{{Z}_{RL}}\left( \omega ,0 \right) \right\} \\
\end{split}
\end{equation}
Taking Laplace transform of Eq.(37) and substituting Eq.(41) into it, we obtain
\begin{equation}
\begin{split}
L\left\{ {}_{0}^{C}D_{t}^{\alpha }f\left( t \right) \right\}&=L\left\{ \int_{0}^{\infty }{{{\mu }_{n-\alpha }}(\omega ){{z}_{C}}\left( \omega ,t \right)d\omega } \right\} \\
 & =\int_{0}^{\infty }{{{\mu }_{n-\alpha }}(\omega ){{Z}_{C}}\left( \omega ,t \right)d\omega } \\
 & ={{s}^{\alpha }}F\left( s \right)+\int_{0}^{\infty }{\frac{{{\left( -\omega  \right)}^{n}}}{s+\omega }{{\mu }_{n-\alpha }}(\omega ){{z}_{RL}}\left( \omega ,0 \right)d\omega }\\
 &  \quad -\int_{0}^{\infty }{{{\mu }_{n-\alpha }}(\omega )\left\{ {{\Delta }_{2}} \right\}d\omega } \\
\end{split}
\end{equation}
where
\begin{equation}
{{\Delta }_{2}}\text{=}\left[ \frac{{{s}^{n-1}}-{{\left( -\omega  \right)}^{n-1}}}{s+\omega } \right]f\left( 0 \right)+\left[ \frac{{{s}^{n-2}}-{{\left( -\omega  \right)}^{n-2}}}{s+\omega } \right]{f}'\left( 0 \right)+\cdots +{{f}^{\left( n-2 \right)}}\left( 0 \right).
\end{equation}

From Eqs. (36)and(43), it is clear that ${{\Delta }_{1}}={{\Delta }_{2}}$, and
\begin{equation}
L\left\{ {}_{0}^{RL}D_{t}^{\alpha }f\left( t \right) \right\}=L\left\{ {}_{0}^{C}D_{t}^{\alpha }f\left( t \right) \right\}.
\end{equation}
Thus,
\begin{equation}
{}_{0}^{RL}D_{t}^{\alpha }f\left( t \right)={}_{0}^{C}D_{t}^{\alpha }f\left( t \right).
\end{equation}
Eq.(45) shows the equivalence of the initialized Riemann-Liouville derivatives and Caputo derivatives with arbitrary orders $\alpha$.

\section{Conclusions}
This paper has presented the equivalence of the initialized Riemann-Liouville derivatives and the initialized Caputo derivatives with arbitrary orders. By synthesizing the initialization function theory and the infinite state theory, the diffusive representations of the two initialized derivatives have been obtained. Laplace transforms of the two initialized derivatives with arbitrary orders have been shown to be equal from Eqs.(18) (30) and (44). As a result, the two most commonly used derivatives have been shown equivalent from Eqs.(19) (31) and (45).\\
\textbf{Acknowledgements}
All the authors acknowledge the valuable suggestions from the peer reviewers. This work was supported by the National Natural Science Foundation of China (Grant No. 11802338).

\end{document}